# DISCUSSION

By Ying Wei[1]

*Columbia University*

First I would like to congratulate the authors for developing a new concept of directional quantile contours. The work will contribute well to the pursuit of multivariate quantiles. The multiple output regression provides a new way of estimating the conditional multivariate quantile functions, which will certainly facilitate a large number of applications. I enjoy the paper for its mathematical rigor and computational savvy. My discussion will focus on the modeling aspect of conditional multivariate quantiles.

**1. Choice of models.** When one tries to incorporate covariate information into multivariate quantiles, certain model assumptions have to be made. As in any regression methods, there are various levels of modeling, for example, linear or nonlinear, parametric or nonparametric. In any application, an appropriate choice of the model matters. I will illustrate this point in my discussion using the same data set as in Hallin, Paindaveine and Šiman (2010). I will later discuss a generalization of the multiple output regression to nonparametric models, and comment on the challenges in model adequacy assessment for the multiple output regression.

To illustrate the main point, let us apply the conditional reference quantiles of Wei (2008) to the same data set in Hallin, Paindaveine and Šiman (2010). The response variables are the calf maximal circumference, denoted as $Y_1$, and the thigh maximal circumference, denoted as $Y_2$. The covariates include age, height, weight and BMI. To make results comparable to those of Hallin, Paindaveine and Šiman (2010), let us estimate the conditional bivariate reference quantile contours of calf and thigh circumferences given the subject's height, weight, age and BMI, separately, as the authors did in their illustrative example. Men and women are analyzed separately. Following the two-step methods of Wei (2008), we first construct stratified

Received August 2009.
[1]Supported by NSF Grant DMS-09-06568 and a career award from NIEHS Center for Environmental Health in Northern Manhattan (ES-009089).

quantile regression models for the conditional joint distribution of calf and thigh circumferences given a chosen covariate $X$. We consider two settings as follows.

1. Setting 1: linear stratified quantile models:
$$Q_{\tau|X}(Y_1) = \alpha_{\tau,1} + \alpha_{\tau,2}X,$$
$$Q_{\tau|X,Y_1}(Y_2) = \beta_{\tau,1} + \beta_{\tau,2}X + \beta_{\tau,3}Y_1.$$

2. Setting 2: nonparametric stratified quantile models:
$$Q_{\tau|X}(Y_1) = \alpha_\tau(X),$$
$$Q_{\tau|X,Y_1}(Y_2) = \beta_{\tau,1}(X) + \beta_{\tau,2}(X)Y_1,$$

where $Q_{\tau|X}(Y)$ denotes the $\tau$th conditional quantile of $Y$ given $X$, and $\alpha_\tau(\cdot)$, $\beta_{\tau,1}(\cdot)$ and $\beta_{\tau,2}(\cdot)$ are smooth functions of $X$, and $\tau \in (0,1)$ is the quantile level.

In both settings, we have a marginal model for $Y_1$, and a conditional model for $Y_2$. Take height, for example, as the covariate, the marginal model in Setting 1 assumes that the $\tau$th quantile of calf maximal circumference $Y_1$ is a linear function of height $X$, while the marginal model in Setting 2 assumes that it is a smooth function of height $X$. Similarly, the conditional model in Setting 1 assumes the quantile of thigh maximal circumference $Y_2$ is linear with both calf maximal circumference $Y_1$ and height $X$, while Setting 2 allows a much more general form. The stratified models in Setting 1 are comparable to the multiple output regression defined in (6.1) of Hallin, Paindaveine and Šiman (2010), since both assume linear structures. More specifically, the stratified models correspond to the directional regression quantiles in Definition 6.1 with $\{\mathbf{b}_{\tau\mathbf{y}} = 0, \mathbf{u}_\mathbf{y} = (1,0)\}$ and $\{\mathbf{b}_{\tau\mathbf{y}} = 0, \mathbf{u}_\mathbf{y} = (1, \beta_{\tau,3})\}$, respectively. In other words, the linear stratified quantile models assume linearity in two specific spatial directions. If we switch the order of $Y_1$ and $Y_2$, we can then obtain another set of stratified models that correspond to another two spatial directions $\mathbf{u}_\mathbf{y} = (0,1)$ and $\mathbf{u}_\mathbf{y} = (\beta_{\tau,3}, 1)$. The two sets of linear models may lead to different approximations. We refer to Wei (2008) for a discussion on the selection and combination of those two possible orders of variables. The multiple output regression makes stronger model assumptions by assuming linearity in all the spatial directions. Because of this stronger model assumption, the multiple output regression is invariant with respect to the ordering of $Y_1$ and $Y_2$.

We fit the data with the stratified models at 200 evenly spaced quantile levels under both settings, and then estimate the 0.2th, 0.5th, 0.8th, 0.94th and 0.98th conditional quantile contours of calf and thigh circumference given the 0.1th, 0.3th, 0.5th, 0.7th, 0.9th quantiles of the covariate, using the model-based simulation approach of Wei (2008). The choice of quantile levels and covariate values match those used in Hallin, Paindaveine and Šiman (2010).

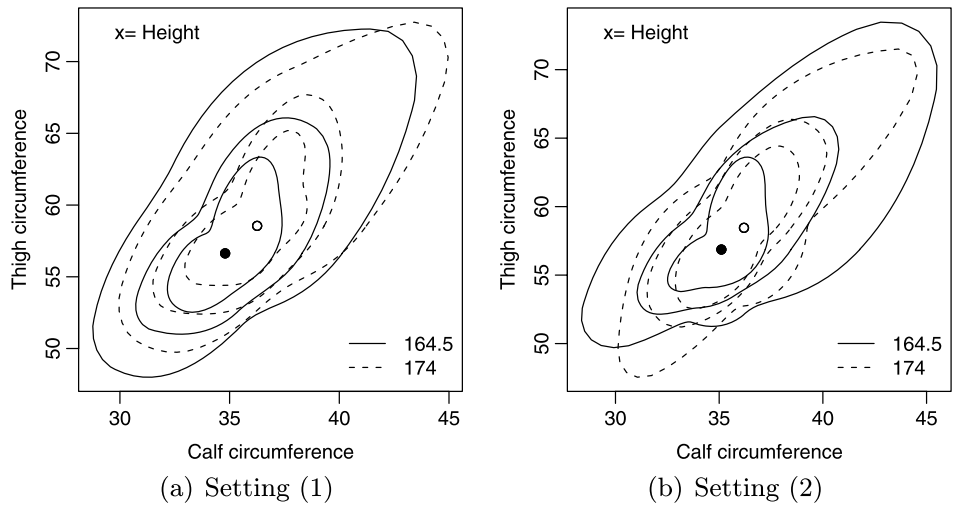

FIG. 1. *Bivariate quantile contours of women's calf and thigh circumferences given different heights. At $X = 164.5$ cm and $X = 174$ cm, quantile contours with $\tau = 0.5, 0.8$ and $0.98$ are shown. The solid contours correspond to $X = 164.5$, and the dotted ones correspond to $X = 174$. Their centers are shown as solid point and open circle, respectively.*

**2. Results.** The resulting reference quantile contours of women's calf and thigh circumferences based on Setting 1 are comparable to Figure 7 of Hallin, Paindaveine and Šiman (2010), but there are noticeable differences in the estimated quantile contours between the linear models and the nonparametric models when height is the covariate.

In Figure 1, we plot the estimated 0.5th, 0.8th and 0.98th reference quantile contours of women's calf and thigh circumferences, conditional on heights at the 0.5th quantile (solid contours) and the 0.9th quantile (dotted contours). Based on linear models (Setting 1), the quantile contours of calf and thigh circumferences of tall women (height = 174 cm) shift upward from those with median height (height = 164.5 cm), which suggests that taller women tend to have larger calf and thigh circumferences than the median-height women. However, that is not true based on the quantile contours generated from the nonparametric models in Setting 2, under which the 0.5th and 0.8th quantile contours of tall and median-height women are fairly close with each other, and the upper part of the 0.98th quantile contours of tall women is actually contained in that of the median-height women. That means that taller women are actually less likely to have really large calf and thigh circumferences relative to the median-height women. This phenomena is even more evident when we analyze men's data. As presented in Figure 2, based on Setting 2, the 0.98th quantile contour of calf and thigh circumferences of taller men (height = 188 cm) are much lower than that of the median height men (height = 177.8 cm). In contrast, based on Setting 1, the distributions of calf and thigh circumferences are comparable for the tall and median-height men. Based on those results, we conjecture that the conditional joint distribution and quantiles of calf and thigh circumference are not linear in height. Consequently, the linear assumptions made in Setting 1, as well as (6.1) in Hallin, Paindaveine and Šiman (2010), lead to biased conclusions.

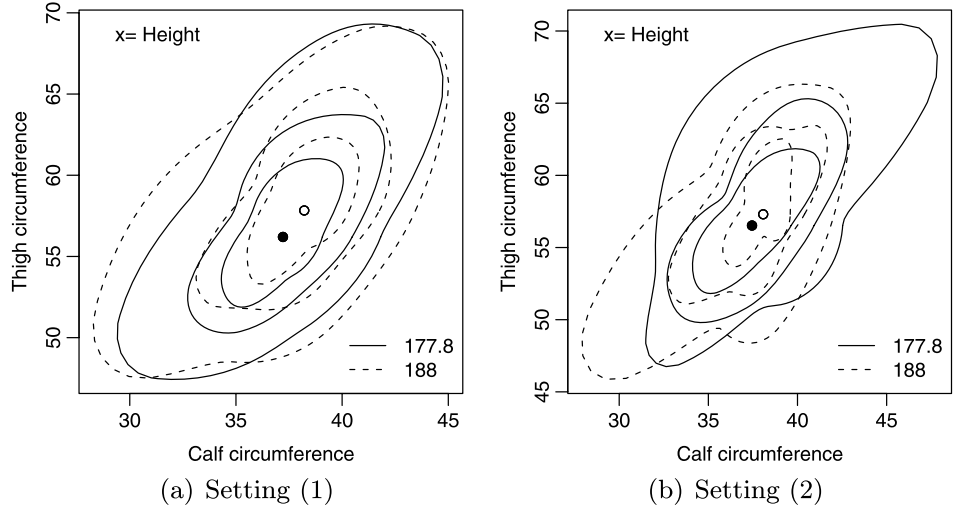

Fig. 2. *Bivariate quantile contours of men's calf and thigh circumferences given different heights. At $X = 177.8$ cm and $X = 188$ cm, quantile contours with $\tau = 0.5, 0.8$ and $0.98$ are shown. The solid contours correspond to $X = 177.8$, and the dotted ones correspond to $X = 188$. Their centers are shown as solid point and open circle, respectively.*

To further support this conjecture, we evaluate the model fitness of the two sets of stratified models, by comparing the model-estimated joint distribution of calf and thigh circumference to the empirical one at the 0.9th quantile of height. The empirical distribution is calculated based on a subsample, consisting of those men whose heights are within a small window of $188 \pm 3$ cm, and the model-based joint distribution is estimated following (2.9) of Wei (2008). The resulting P–P plot is presented in Figure 3. The P–P plot depicts how close the model estimated joint distribution is to the empirical one. If the models fit the data well, the two distributions should be close to each other. Based on Figure 3, it is clear that, conditional on height 188 cm, the model-estimated distribution under Setting 1 over-estimated the upper quantiles of calf and thigh circumferences, which in turn indicates the lack-of-fit of the linear models. In the mean time, the nonparametric models provide a good approximation to the conditional joint distributions given height 188 cm, as shown in Figure 3(b). Based on the discussion above, we believe that the linear models are not adequate for the conditional quantiles of calf and thigh circumferences given height.

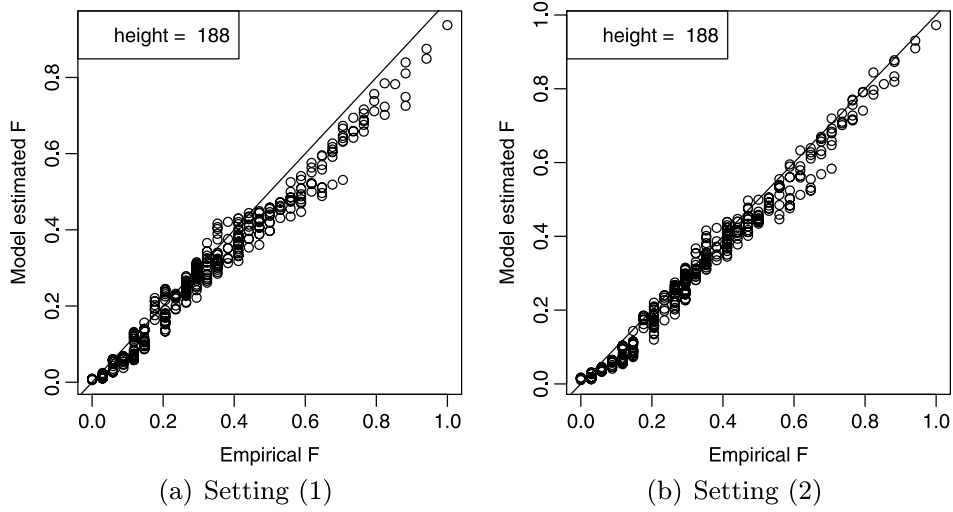

(a) Setting (1)     (b) Setting (2)

Fig. 3. *PP-plot for assessing the conditional model fitness for men's calf and thigh circumferences given height* $= 188$ *cm.*

**3. Nonparametric multiple output regression.** Similar to the linear models in Setting 1, the definition of multiple output regression in (6.1) also assumes linearity between the response $Y$ and the covariate $X$. Consequently, as shown in the previous example, it may not be adequate to model the conditional quantile contours in some applications. Therefore, it might worth the efforts to extend the linear multiple output regression to more general cases. If there is only a single covariate as illustrated in the example, one natural extension is to replace the covariate $X$ by its B-spline basis functions. That is, one can replace definition (6.1) by

$$(1) \qquad (\mathbf{b}_{\tau \mathbf{y}}, g_\tau(x)) = \arg \min_{\mathbf{b_y}, g(x)} E[\rho_\tau(\mathbf{u}'_\mathbf{y} \mathbf{Y} - \mathbf{b}'_\mathbf{y} \mathbf{\Gamma}'_\mathbf{u} \mathbf{Y} - g(X))],$$

where $g_\tau(X)$ is an unknown smooth function of the covariate $X$. The function $g_\tau(X)$ is to be approximated by $g_\tau(x) \approx \mathbf{b}'_{\tau \mathbf{x}} \pi(X)$, with $\pi(X)$ being $q$-dimensional B-spline basis functions given certain internal knots and order of spline. The solution is then in the following form,

$$(2) \qquad (\mathbf{b}'_{\tau \mathbf{y}}, \mathbf{b}'_{\tau \mathbf{x}})' = \arg \min_{\mathbf{b_y}, \mathbf{b_x}} E[\rho_\tau(\mathbf{u}'_\mathbf{y} \mathbf{Y} - \mathbf{b}'_\mathbf{y} \mathbf{\Gamma}'_\mathbf{u} \mathbf{Y} - \mathbf{b}'_\mathbf{x} \pi(X))],$$

and the outlined linear programming algorithm in Hallin, Paindaveine and Šiman (2010) can be applied directly. If there are more than one covariates, then an additive model can be considered.

**4. Model assessment of multiple output regression.** As various models of different complexity can be considered for the multiple output regression, it is helpful to evaluate model adequacy. Typical goodness-of-fit statistics for multivariate quantile contours may not be applicable to the quantile contours using the multiple output regression, since the directional quantile contours do not have the coverage property in the sense of Serfling (2002). Unlike other definitions of multivariate quantile functions, in which the $\tau$th quantile contour usually has the coverage probability $\tau$, the probability mass of the $\tau$th directional quantile contour, however, is actually unknown. For the same reason, they do not have a direct mapping to the distribution functions. However, since those directional quantile contours are generated from the regression quantiles at all the spatial direction $\mathbf{u_y}$, one may assess the model adequacy of the specified models over all the spatial directions. Take the linear model, for example, and suppose that $(\mathbf{x}_i, \mathbf{y}_i), i = 1, \ldots, m$, is a subset where all the $\mathbf{x}_i$ are equal to or close to a target value of $x$. If the specified model fit the data well at the direction $\mathbf{u}$ given the covariate $x$, then we expect that

$$\Delta(\mathbf{u}, x) = \frac{\sum_{i=1}^m I\{\mathbf{u}'_\mathbf{y} \mathbf{y}_i - \widehat{a} - \widehat{\mathbf{b}}'_y \mathbf{\Gamma}'_\mathbf{u} \mathbf{Y} - \widehat{\mathbf{b}}'_\mathbf{x} x_i \leq 0\}}{m} \approx \tau,$$

where $(\widehat{a}, \widehat{\mathbf{b}}_\mathbf{y}, \widehat{\mathbf{b}}_\mathbf{x})$ are the estimated coefficient. An overall model adequacy measure can then be constructed by integrating $\Delta(\mathbf{u}, x)$ over the entire spatial directions, that is,

$$\Delta(x) = \int_u \{\Delta(\mathbf{u}, x) - \tau\} d\mathbf{u}.$$

If the multiple output regression holds, then $\Delta(x)$ should be close to zero. Further research is clearly needed to make this diagnostic tool broadly useful.

**5. A final note.** Hallin, Paindaveine and Šiman (2010) compared the directional quantile contour with the reference quantile contour of Wei (2008). The authors are right in pointing out that the reference quantile contours would depend strongly on the choice of the center, but it is also worth noting that Wei (2008) uses the component-wise medians as a specific choice of centers for the reference quantile contours. This way, the definition of reference quantile contours reduced to the reference percentile charts [Cole and Green (1992)] for one dimensional $Y$. Because no single approach to multivariate quantile contours is likely to be the best in all applications, the ideas proposed by Hallin, Paindaveine and Šiman (2010) is an exciting addition to an important area of multivariate quantiles.

Department of Biostatistics
Columbia University
New York, New York 10032
USA
E-mail: ying.wei@columbia.edu